\documentclass[11 pt]{article}
\usepackage{mathptmx}%
\usepackage{indentfirst} 

\usepackage{titling}
\preauthor{\begin{center}
		\large \lineskip 0.5em%
				\vspace{-3.1cm}
						\begin{tabular} [t] {c} } 
	\postauthor{\end{tabular} \par \end{center}}

\usepackage{caption}

\usepackage{longtable}
\usepackage{graphicx}
\usepackage{url}
\usepackage{amsmath}
\usepackage{lscape}
\usepackage{mathrsfs}
\usepackage{subcaption}
\usepackage{titlesec}
\titlelabel{\thetitle.\;}
\titleformat*{\section}{\large\bfseries\sffamily}
\titleformat*{\subsection}{\large\itshape}

\usepackage[labelsep=period]{caption}
\usepackage[paperwidth=170mm,paperheight=240mm]{geometry}
\author{\and Frank Lad\thanks{University of Canterbury, Department of Mathematics and Statistics, Christchurch, New Zealand, {\it E-mail address:} F.Lad@math.canterbury.ac.nz}  \and  - \and  Giuseppe Sanfilippo\thanks{Department of Mathematics and Computer Science, University of Palermo, Italy, {\it E-mail address:} giuseppe.sanfilippo@unipa.it}}

\title{\vspace{-1cm}  \large  \bf SCORING ALTERNATIVE FORECAST DISTRIBUTIONS: \ COMPLETING THE KULLBACK DISTANCE COMPLEX }

 \date{   }

\begin{document}
	\clearpage\maketitle
	\thispagestyle{empty}
\pagestyle{empty}
\section*{Abstract} 
We develop two surprising new results regarding the use of proper scoring rules for evaluating the predictive quality of two alternative sequential forecast distributions. Both of the proponents prefer to be awarded a score derived from the other's distribution rather than a score awarded on the basis of their own. 
A Pareto optimal exchange of their scoring outcomes provides the basis for a  comparison of forecast quality that is preferred by both  forecasters, and also evades a feature of arbitrariness inherent in using the forecasters' own achieved scores.  The  well-known Kullback  divergence, used as a measure of information, is evaluated via the entropies in the two forecast distributions and the two cross-entropies between them.  We show that 
Kullback's symmetric measure needs to be appended by three component measures if it is to characterise completely the information content of the two asserted probability forecasts.  Two of these do not involve entropies at all.  The resulting ``Kullback complex'' supported by the 4-dimensional measure is isomorphic to an equivalent vector measure generated by the forecasters' expectations of their scores, each for one's own score and for the other's score. 
We foreshadow the results of a sophisticated application of the Pareto relative scoring procedure for actual sequentional observations, and we propose a standard format for evaluation.\\


\noindent {\bf JEL Classification:} C10, 
C11\\

\noindent {\bf Key Words:} Total logarithmic scoring rule; prevision;  entropy/extropy;  cross entropy; Pareto optimal exchange; Kullback symmetric divergence;  Bregman divergence. 


\section{Introduction} 
Bruno de Finetti's ``partly-baked'' musings with I. J. Good (1962) on the use of a penalty function to gauge the relative predictive 
qualities of competing forecast distributions have developed by now into  a complete theory of proper scoring rules, replete with practical computational procedures.  Gneiting and Raftery (2007) present a fairly comprehensive technical review of both the theory and its application. Lad (1996, Chapter 6) presents a pedagogical introduction.  Proper scoring procedures are meant to provide the basis for a subjectivist understanding of forecast evaluations that can replace completely the common practice of ``hypothesis testing.''  Such testing is regarded by subjectivists to be without foundation.  The very notion that data are generated by an unobservable probability distribution is rejected in their characterisation of probability and its connection with  empirical statistical research.  In this present tribute to the imagination of de Finetti we propose a striking development in scoring theory.  Our analysis gives rise to a completed understanding of Kullback's information measure in the context of the total logl scoring function.  Implications for statistical practice are substantial.\\


The logarithmic scoring function for a probability mass vector (pmv) ${\bf p}_N$ asserted by a forecaster ``p'' for a quantity $X$ on the basis of an observation $X = x^o$ within its realm of possibility $\mathscr{R}(X) = \{x_1, x_2, ..., x_N\}$  is specified as the logarithm of the probability assessed for the outcome $x^o$ that is actually observed:  \ $S(X = x^o, {\bf p}_N) \equiv log(p^o)$. \ It has been known since early developments of proper scoring rules by Savage (1971) and by Matheson and Winkler (1976) that a pmv proponent's  expectation of one's own log score equals the negentropy value embedded in that pmv.
Symbolically,  $E_p[log\;P(X=x^o))] \;  = \; \sum p_i \;log\: p_i$.  Moreover, p's expectation of the amount by which p's own log score will exceed the score of someone else's assertions, ${\bf q}_N$, is equal to the relative entropy in ${\bf p}_N$ with respect to ${\bf q}_N$.  That is, 
\[E_p[S(X,{\bf p}_N) - S(X,{\bf q}_N)] \; = \; \sum p_i\;log\:(\frac{p_i}{q_i}).\]
This relative entropy measurement is also widely known as the Kullback-Leibler (1951) directed divergence between the pmv's ${\bf p}_N$ and ${\bf q}_N$, and is commonly denoted by $D({\bf p}_N||{\bf q}_N)$.  A non-symmetric function of the two pmv's, the value of $D({\bf p}_N||{\bf q}_N)$ is summed with its counterpart $D({\bf q}_N||{\bf p}_N)$ to yield Kullback's (1959) symmetric distance measure, denoted by $\mathscr{D}(p||q) \: \equiv \: D({\bf p}_N||{\bf q}_N)\,  + \,D({\bf q}_N||{\bf p}_N)$.  The analysis and the applications of this distance function in information theory and practice have been extensive.  Among recent theoretical developments, the Kullback-Leibler function is now known to constitute a ``Bregman divergence'' relative to the Bregman function $\Phi({\bf p}_N) = \sum p_i \;log\: p_i$, studied by Censor and Zenios (1997).  This structure sets it in a context shared by several other divergence functions, though not by still others.  The recent contribution of Dawid and Musio (2015) references how lively and widespread is research interest in related issues.\\

Subsequent consideration of scoring rules by Lad, Sanfilippo and Agr\`o (2012) 
has addressed the fact  that the log score for mass functions, while proper, is incomplete in that it does not make use of {\it all} the observations with their asserted probabilities that are entailed in the assertion of ${\bf p}_N$ and the observation of $X$. 
Augmented by the sum of $(1-p_i)\;log(1-p_i)$ values for possibilities of $X$ that {\it are} observed {\it not} to obtain, the log scoring rule becomes the {\it total log scoring rule}.  Its prevision (expectation) according to the proponent who asserts ${\bf p}_N$ equals the negentropy {\it plus} the negextropy of the asserted mass function.  Negextropy was characterized in an  article by Lad, Sanfilippo, and Agr\`o (2015) as a complementary dual of negentropy.  The structure of the Kullback-Leibler divergence entertains a complementary dual as well: \ $D^c({\bf p}_N||{\bf q}_N) = \sum (1-p_i) log(\frac{1-p_i}{1-q_i}). $
Parallel considerations define the {\it complementary} symmetric Kullback distance as \; $\mathscr{D}^c(p||q) \equiv D^c(p||q) + D^c(q||p)$.  Together, the Kullback distance and its complementary distance sum to the Kullback total symmetric distance: $\mathscr{D}_T(p||q) \; \equiv \; \mathscr{D}({\bf p}_N||{\bf q}_N) \; + \;\mathscr{D}^c({\bf q}_N||{\bf p}_N)$. \\

The present article presents two new developments in our understanding of associated issues.  In the first place, we consider the expectations of the asserters of both ${\bf p}_N$ and ${\bf q}_N$, each for one's own total log score {\it and} for the other's score.  
Attention turns to these proponents' considerations of the net gains they each might achieve on the basis of their professed willingness to trade their achieved scores for their personally assessed prevision values. (This constitutes the very definition of their previsions.) We find an appropriate strategy for assessing the relative merits of the two pmv assertions: \   
perhaps surprisingly, it is to award to each forecast proponent the {\it net gain achieved by the other's} pmv, scaled appropriately, rather than to award to each forecaster one's own achieved proper score. \\

Secondly, our analysis turns to the pair of expected net gains and expected comparative gains of the forecasters.  We learn that Kullback's total symmetric distance measure can be generated equivalently by two other pairs of functions that are different from that which specifies the Kullback-Leibler divergence.  This result suggests that Kullback's symmetric distance function should be extended to a four-dimensional vector function if it is to characterise all the distinct features of the information expressed in the two pmv's.  It turns out that the vector space spanned by this 4-D vector function is isomorphic to the space spanned by the two proponents' previsions for their own scores and for each other's score.\\



Section 2 of this article identifies the two pmv proponents' expected {\it net gains} and expected {\it comparative gains} according to their uncertain assessments of an unknown quantity $X$. Section 3 develops the logic of the Pareto-optimal assessment of their comparative forecast quality we propose here, while the {\it scale} of this comparative scoring procedure is discussed  in Section 4. Scaling requires consideration of the variances assessed by the two forecasters for various score components.   Attention turns to the higher dimensions of the Kullback complex that arise from these considerations in Section 5, and Section 6 displays why the appended dimensions of the complex are {\it not} Bregman divergences.   Section 7 portrays the historical context of Kullback's analysis of the problem, and outlines its reconsideration according to the subjectivist understanding.  Section 8 presents briefly some  experience with and some proposals for application to statistical problems.\\

Deliberations in this article make use of the syntax preferred for Bruno de Finetti's operational subjective construction of coherent prevision, discussed by de Finetti (1967, 1972) and Lad (1996, pp. 39-43).  Events are characterised as numerical quantities (numbers, not sets) with the realm of possibility $\mathscr{R}(E) = \{0, 1\}$.  Otherwise they are addressed with the same syntax as any other quantity whose realm is finite and discrete.  Prevision, $P(X)$, is a personal value assertion about $X$, a price at which its proponent is avowedly indifferent to a prize of $X$ or of $P(X)$.  Since a prevision is an assertion made by someone, the $P$ symbol shall be subscripted by a mark denoting who is making the assertion, as in $P_p(^.)$ and $P_q(^.)$ used here.  We shall refer to the proponents of these previsions merely as ``p'' and ``q'', respectively.  One's prevision for an event is also called a probability.  Parentheses around an expression that may be true or false denotes an event whose numerical value equals $1$ if the expression is true, and $0$ if it is false.  Thus, for example, $E_i \equiv (X=x_i)$ is an event that equals $1$ if and only if it happens that $X$ is in fact found to equal $x_i$;  otherwise $E_i = 0$.  In the more commonly used measure-theoretic syntax of formalist probability ala Kolmogorov, prevision is equivalent to expectation;  in these terms, the expectation of a $0-1$ random variable is a probability.

\section{Interpersonal assessment of the Total Log Scores}

The purpose of this Section is to formalise the notation for the total logarithmic score that will be used to assess our two forecasting distributions (pmvs), and to introduce the formalities of two concepts that will be central to our considerations.  The first is the {\it net} gain achieved by a forecaster who asserts a pmv for an observable quantity.  It defined as the difference between the score that will be achieved and the forecaster's expectation for the score.  The second is the {\it comparative} gain achieved by one of the forecaster's score relative to that achieved by the other.  It is defined simply as the difference between their two scores.  We can then identify each forecaster's professed expectation for each of these gains.  In doing so we will have identified {\it four fundamental previsions} that are relevant to evaluating the forecast quality of the two forecast pmv's on the basis of observations.\\  

Consider the standard problem of statistical forecasting an unknown quantity $X$ whose realm of possible observation values is $\mathscr{R}(X) = \{x_1, x_2, ..., x_N\}$.  An associated vector of partition events, {\bf E}$_N$, is composed of its logically related components $E_i \equiv (X=x_i)$ for $i = 1, 2, ..., N$.  
The sum of these events must equal $1$ since one and only one of them must equal $1$.  Because the considered value of $X$ is unknown when the forecasters assess probabilities, it is just not known which of the events is the $1$ and which of them are the $0$'s.  This will be identified when we observe the value of $X$. 
We shall study the uncertain opinions of two people who publicly assert distinct probability mass functions via the vectors ${\bf p}_N$ and ${\bf q}_N$, respectively.  Once the value of $X$ is observed, their assertions will be assessed via a widely touted proper score, the Total Log Scoring function, defined for {\bf p}$_N$ (and similarly for {\bf q}$_N$) by 
\begin{eqnarray}
S_{TL}(X, {\bf p}_N) \ &\equiv& \ \sum_{i=1}^N E_i\;log(p_i) + \sum_{i=1}^N {\tilde E}_i\;log(1-p_i) \nonumber  \\ 
&=& \ log(p^o) + \sum_{i\,:\,p_i\neq p^o} log(1-p_i)
\nonumber\\
&=&log(\frac{p^o}{1-p^o}) \; + \;\sum_{i=1}^N \;log(1-p_i) \ ,
\label{eq:totlogscore}\end{eqnarray}
\noindent as long as the pmv {\bf p}$_N$ is strictly positive and $N\;>\;1$.  Here $p^o$ denotes the probability asserted for the value of
$(X=x^o)$ that eventually comes to be observed; $\ $and the sum in the second line of $(\ref{eq:totlogscore})$ includes all summands of the form 
$log(1-p_i)$ {\it except for}\: $log(1-p^o)$.  \\

The ``total log score'' is a
negative valued score, but it constitutes an award.  A larger score (smaller absolute
value) is achieved by a more informative distribution.  As extremes, the pmv that
accords probability $1$ to the value of $X$ that occurs, $(X=x^o)$,
receives a score of
$0$; \ those which accord probability $1$ to any particular value of $X$ that does {\it not} occur receive $-\infty$. \ The uniform pmv receives $log(1/N) + (N-1)\;log(1-1/N) \approx -log(N) + e^{-1}$ for large $N$.
\vspace*{.33cm}

\noindent {\bf N.B.} \ From this point on we shall no longer write subscripts on vector quantities to denote their size, nor the limits on the range of summations, always understood to be summations over units from $1$ to $N$.\\



To foreshadow a relevant consideration that will be addressed at the very end of the analysis, we note here that each of the two summands constituting the total log score, $\sum p\;log\,p$ and $\sum (1-p)\;log(1-p)$ is itself a proper score of the pmv ${\bf p}$, as is any positive linear combination of these two.  Whereas the constructions we now discuss are developed in terms of the total log score, they pertain to each of these addends separately in the same way.

\subsection{Previsions for one's own score and for another's score}

In the following discussion we shall denote the entropy and extropy functions, respectively, by $H({\bf p}) = -\sum p\;log(p)$ and $J({\bf p}) = -\sum (1-p)\;log(1-p)$.  Sums such as $-\sum p\;log(q)$ and $-\sum q\;log(p)$ are called ``cross entropies'', and will be denoted by $CH({\bf p,q})$ and $CH({\bf q,p})$.  Similarly, ``cross extropies'' will be denoted by $CJ({\bf p,q}) \equiv -\sum (1-p)\; log(1-q)$.  Recognise that cross entropies/extropies are not symmetric in their arguments.  For example, $CH({\bf p,q}) \neq CH({\bf q,p})$. \\

\noindent {\bf Theorem 1.} \ \ The linearity of coherent prevision identifies the previsions of both p and q for their own scores and for the other's score as:
\begin{small}
\begin{eqnarray}
P_p[S_{TL}(X,{\bf p})] \, &=& \, \sum p\; log(p)  \ + \ \sum (1-p)\; log(1-p)  =  -[H({\bf p})+J({\bf p})] \ , \ {\rm  and } \nonumber\\
P_p[S_{TL}(X,{\bf q})] \, &=& \, \sum p\; log(q)  \ + \ \sum (1-p) \; log(1-q) =  -[CH({\bf p,q}) + CJ({\bf p,q})]\; , \nonumber \\
P_q[S_{TL}(X,{\bf q})] \, &=& \, \sum q\;  log(q) \ + \ \sum (1-q) \; log(1-q)  =  -[H({\bf q})+J({\bf q})]\; ,\  {\rm \ and} \nonumber \\
P_q[S_{TL}(X,{\bf p})] \,  &=& \, \sum q\;log(p)  \ + \ \sum (1-q) \; log(1-p)  =   -[CH({\bf q,p}) + CJ({\bf q,p})]\; .\nonumber \label{eq:4fundprevs}
\end{eqnarray}
\end{small}
\noindent The proof is immediate from the application of linear coherent previsions to the first line of equation $(\ref{eq:totlogscore})$. \vspace{.22cm}

These four previsions exhaust the content of what p and q assert about their own prospective score values and for the scores to be achieved by one another, before the value of $X$ is observed.  We shall designate them as {\it the four fundamental previsions}, and we shall learn why this nomenclature is appropriate.  Let's think a moment about their relative sizes.\\

Since the scoring function $S_{TL}(^.,^.)$ is proper,  each of p and q expects to achieve a greater score than will the other;  for the proponent of a pmv expects to achieve a maximum score by publicly professing one's actual pmv, {\bf p}, rather than any other pmv, {\bf q}.  For this reason, the use of proper scoring rules is said to promote honesty and accuracy in the profession of probabilities.  Thus, for example, $-[H({\bf p})+J({\bf p})]$ surely exceeds $-[CH({\bf p,q}) + CJ({\bf p,q})]$, a statement which is true when {\bf p} and {\bf q} exchange places as well.  However, the relative sizes of the two assessors' own expected scores, $P_p[S_{TL}(X,{\bf p})]$ and $P_q[S_{TL}(X,{\bf q})]$, may be in any order.  This just depends on what they each assert about $X$, which is {\bf p} or {\bf q}, respectively. \\

When cross-score previsions are considered, p's prevision for q's score may be {\it either smaller or larger} than q's expectation for q's own score.  In a computational survey of pmv's within the unit-simplex we found that the size of $P_p[S(X,q)]$ exceeds $P_q[S(X,q)]$  in about 10 to 15 percent of paired $({\bf p}_3, {\bf q}_3)$ selected uniformly at random from the unit simplex.  We have not yet characterised the situations in which this occurs, though it should be straightforward.
\vspace*{.1cm}

\noindent {\bf N.B.} \ From this point we no longer use bold print to distinguish vectors from their  components.  These will be identifiable by their context. \vspace*{.3cm}

We are now ready to analyse important features of p's and q's expectations regarding one another's scoring performance.

\subsection{Expected Net Gains and Expected  Comparative Scores}

Consider first the ``net gains'' to be achieved by the proponents of p
and q as a result of their receiving the total log scores, $S(X,p)$ or $S(X,q)$ respectively, as a gain in return for paying out the prices at which each of them values this score: \begin{eqnarray} NG(X,p) \ &\equiv& \  S(X,p) - P_p[S(X,p)]\ \ , \ \ {\rm and} \nonumber\\
NG(X,q) \ &\equiv& \  S(X,q) - P_q[S(X,q)]\ . \ \ \ \ \label{eq:netgains}
\end{eqnarray}
Both $p$ and $q$ expect a personal net gain of $0$, since their assertions are presumed to be
coherent.  One's prevised (expected) score is the price at which a person values this unknown score to be achieved.  The proponent of a pmv is avowedly willing both to buy and to sell a claim to the score for this price.\  For example, p is avowedly indifferent to the values of $S(X,p)$ and $P_p[S(X,p)]$. Thus, \ $P_p[NG(X,p)] = 0$, as does  $P_q[NG(X,q)] = 0$.  \\



However, the proponents of p and q do {\it not} expect each others' net gains to equal $0$.  Before assessing them, it is worthwhile to recognize that there is something arbitrary about our definition of Net Gain.  For the forecaster is also indifferent to an award of the negative net gain, the result of an exchange in which the forecaster offers to sell the value of $S(X,p)$ in return for receiving $P_p[S(X,p)]$.  This would reverse the sign of an awarded Net Gain as we have defined it.  This is the sense in which an award to p of a proper score value or its prevision assessment would be arbitrary.   Nonetheless, we are now ready to study it as defined.\\

\noindent {\bf Theorem 2.}  $P_p[NG(X,q)]$ may be positive or negative valued, as may $P_q[NG(X,p)]$.\\

\noindent As proof, note specifically that 
\begin{eqnarray} P_p[NG(X,q)]  &=& P_p[S(X,q)] - P_q[S(X,q)] \nonumber =\\
&=&  -[CH(p,q) + CJ(p,q)] + H(q) + J(q) \nonumber \\
&=& \sum (p-q)\; log(\frac{q}{1-q})\; ; \nonumber \\
{\rm while}\vspace{.35cm}\ \ \ \ \ \ \ \ \ & & \label{eq:prevotherNG}\\
\noindent P_q[NG(X,p)] &=& P_q[S(X,p)] - P_p[S(X,p)]  =\nonumber \\
&=&  -[CH(q,p) + CJ(q,p)] + H(p) + J(p)\; \nonumber\\
&=& \sum (q-p)\; log(\frac{p}{1-p})\; \ \ . \nonumber 
\end{eqnarray}

\noindent Now we have already noticed and remarked that p's prevision for 
q's score may be greater than or less than q's prevision for the same.  If p expects a score for q greater than q expects, 
then p will expect a positive net gain for q because p is thinking that q expects too low a score.  In such a case, 
\[P_p[NG(X,q)] = P_p[S(X,q)]-P_q[S(X,q)] \ \  > \ \  P_q[S(X,q)]-P_q[S(X,q)] = 0.\] 
 If p expects a score for q less than q expects, the reverse is be true: \ $P_p[NG(X,q)] = P_p[S(X,q)]-P_q[S(X,q) ]\  < \ 0.$\\
 

The forecasters' expectations of their {\it comparative} gains relative to the other's are another matter.  Defining the ``comparative gain of p over q'' on the basis of observing $X$ as \begin{eqnarray}
CG(X,p,q) \  &\equiv& \  S(X,p)- S(X,q) \; = \; -CG(X,q,p) \ \ , \label{eq:compgain}
\end{eqnarray} it is easy to see that both forecasters expect a positive comparative gain.\\

\noindent {\bf Theorem 3.}  \ Both forecasters always expect a positive comparative gain for themselves.\\

\noindent As proof, it is clear that
\begin{eqnarray} P_p[CG(X,p,q)] \, &=& \, P_p[S(X,p)] - P_p[S(X,q)]= \,\nonumber \\
 &=&  \,-H(p) + CH(p,q) - J(p) + CJ(p,q)
\, ; \nonumber\\
{\rm and \ similarly,} & & \label{eq:prevscompgain}\\
P_q[CG(X,q,p)] \, &=& \, P_q[S(X,q)] - P_q[S(X,p)]\, =\nonumber \\
&=& \, -H(q)+ CH(q,p) - J(q)  + CJ(q,p)
\, ,\nonumber 
\end{eqnarray}
\noindent using the four fundamental previsions we identified in Theorem 1.  Both of these previsions are positive because the scoring rule is proper.\\

Now comparing the equation pairs numbered $(\ref{eq:prevotherNG})$ and $(\ref{eq:prevscompgain})$, it is apparent that {\it summing} each of the paired expected net gains and comparative gains yields precisely the negative results of one another.  It is worth stating this too as a theorem. \vspace{.22cm}\\
\noindent{\bf Theorem 4.} \ $P_p[NG(X,q)] \; + \; P_q[NG(X,p)] \ =  \ \vspace{.1cm}\\
\hspace*{3cm}- \; \{P_p[CG(X,p,q)] \; + \; P_q[CG(X,q,p)]\}\vspace{.1cm}\ .$ \\
\noindent Notice however that the paired summands for these opposite magnitudes are different.\vspace{.4cm}

We shall now continue this analysis by devising an improved procedure for the method of comparative scoring, and then we shall use this result to extend our understanding of Kullback's symmetric divergence.

\section{Evaluating pmv's by trading inhering Net Gains}



Since p assesses the personal net gain $NG(X,p)$ with prevision 0 while valuing $NG(X,q)$ with the non-zero prevision specified in equation $(\ref{eq:prevotherNG})$, p would surely be willing to trade a claim to plus or minus $NG(X,p)$ in exchange for a claim to $NG(X,q)$ if $P_p[NG(X,q)]$ is positive, and to $-NG(X,q)$ if that prevision is negative.  For p assesses $NG(X,p)$ with value $0$, whereas either $NG(X,q)$ or $-NG(X,q)$ is assessed with positive value.  In the same way, q values claims to $NG(X,q)$ with prevision $0$, while valuing the net gain of p with the  prevision $P_q[NG(X,p)]$, as specified in in the lower half of equation $(\ref{eq:prevotherNG})$.  Thus both p and q would be pleased to make an exchange in which p receives (appropriately either the positive or negative) value of the $NG(X,q)$ from q while q receives similarly the value of the net gain $NG(X,p)$ (or its negative) from p. In this exchange, both p and q would each be providing the other with something personally regarded as worthless.   Thus, both of them would be pleased by the exchanges.  Altogether they constitute a ``Pareto optimal'' exchange. \\

To clarify this waffle about receiving ``a net gain or its negative'', note for example that q is avowedly willing to buy {\it or to sell} the value of $S(X,q)$ for $P_q[S(X,q)]$, since q values these equally.  Now with these two options available from q as proclaimed, p would willingly take up q's offer to to sell $S(X,q)$ at this price if $P_p[S(X,q)] > P_q[S(X,q)]$, {\it or} to take up q's offer to buy $S(X,q)$ at this price if $P_p[S(X,q)] < P_q[S(X,q)]$.  The same structure of relative preference applies to q in reacting to p's proclaimed indifference to $S(X,p)$ and $P_p[S(X,p)]$.  To simplify the formalities of this waffling about the action whether to buy or to sell, during further discussion we shall denote the content of such exchanges by saying $NG^*(X,q)$ is exchanged for $NG^*(X,p)$.  This may mean that p gives $NG(X,p)$ to q in return for $P_p[NG(X,p)]$ {\it or} that p gives $P_p[NG(X,p)]$ to q in exchange for $NG(X,p)$, at the choice of q.  Simultaneously, q gives to p the value of $NG(X,q)$ in return for $P_q[NG(X,q)]$ or vice versa at the choice of p.  Of course this specificity would need to be recognized and recorded in the computation of the resulting net gains.\\

One difference in the assessed values of the two sides of such an exchange is that both p and q would assess their own net gains being given up with a different variance than the net gains they would receive.  This constitutes no problem for the usual characterisation of coherent prevision nor for the acceptability of the proposed exchanges.  In defining prevision, the scales of exchange are overtly kept small enough that linearity of utility applies to practical valuations.  However, we can rescale these agreeable exchanges so that both p and q would regard the two sides of their exchanges symmetrically with the same variance, that is, the same risk.  To facilitate this redesign, the pmv proponents would want to rescale the two sides of the exchanges they offer by their assessed standard deviations for them.  That is,

\vspace*{.3cm}
\indent p would offer to give q \ \ $\frac{NG^*(x,p)}{SD_p[NG(X,p)]}$ \ \ in exchange for \ \ $\frac{NG^*(x,q)}{SD_p[NG(X,q)]}$ \ \ ;

\vspace*{.3cm}
\noindent while at the same time 

\vspace*{.3cm}
\indent q would offer to give p \ \ $\frac{NG^*(x,q)}{SD_q[NG(X,q)]}$ \ \ in exchange for \ \ $\frac{NG^*(x,p)}{SD_q[NG(X,p)]}$ \ \ .

\vspace*{.3cm}
\noindent Both p and q regard both of these offers with positive value, so both would be accepted.  

\vspace*{.3cm}
Thus, we have found an intriguing and novel solution to the question of how to evaluate the relative quality of p's and q's pmv assertions when the value of $X$ is actually observed.  Simplifying the net result of those two exchanges that are agreeable to both of them, the appropriate way to score the pmv assertions of p and q would be 

\vspace*{.3cm}
\indent to award to q the score of \ \ $\frac{NG^*(x,p)}{SD_p[NG(X,p)]} \ + \ \frac{NG^*(x,p)}{SD_q[NG(X,p)]}$, \ \ \ and \vspace*{.3cm}\\
\indent to award to p the score of \ \ $\frac{NG^*(x,q)}{SD_p[NG(X,q)]} \ + \ \frac{NG^*(x,q)}{SD_q[NG(X,q)]}$. \\

\noindent These exchanges constitute  for us a {\it Pareto optimal scoring procedure for  comparing the relative quality of p's and q's forecasting distributions}\,:\ \ to award to q the amount of p's achieved Net Gain$^*$ scaled by the sum of the square roots of p's and q's precisions for that Net Gain$^*$, and to award to p the amount of q's Net Gain$^*$ scaled by the sum of the square roots of p's and q's precisions for {\it that} Net Gain$^*$.\\


In order to compute these relative information scores, we need to compute four relative variances:  $V_p\{NG[S(X,p)]\}, V_p\{NG[S(X,q)]\}, V_q\{NG[S(X,p)]\}$, and $V_q\{NG[S(X,q)]\}$.  Let us  digress momentarily to think about these variances in a brief Section of its own.

\section{... and as to the variances}

Here are some considerations regarding the cross assessments of the standard deviations.  Firstly, remember that $NG(X,p) = S(X,p) - P_p[S(X,p)]$.  So $V_p[NG(X,p)] = V_p[S(X,p)]$, because  $P_p[S(X,p)]$ is a specified number.  Now think of the total log score in the form of $S(X,p) \; = \; log(\frac{p^o}{1-p^o}) \ + \ \sum log(1-p)$ identified in the third line of equation $(\ref{eq:totlogscore})$.  In this form it is clear that \vspace{.3cm}\\
\indent $V_pS(X,p) = V_p[log(\frac{p^o}{1-p^o})]$\; , \ \ and similarly, \ \
$V_pS(X,q) = V_p[log(\frac{q^o}{1-q^o})]$\; .\vspace{.3cm}\\
\noindent These follow for the same reason, that the summations $\sum log(1-p)$ and $\sum log(1-q)$ are specified numbers.\\

\noindent Computationally, \ \ \ \ \ $V_p[S(X,p)] \ = \ \sum\;p_i\;[log(\frac{p_i}{1-p_i})]^2 \ - \ [\sum\;p_i\;log(\frac{p_i}{1-p_i})]^2\;,$ \\

 while \ \ \ \ \ \ \  \ \ \ \ \ \ \ \ \; $V_p[S(X,q)] \ = \ \sum\;p_i\;[log(\frac{q_i}{1-q_i})]^2 \ - \ [\sum\;p_i\;log(\frac{q_i}{1-q_i})]^2\;$ \vspace{.3cm} \\
 \noindent according to the standard result for any quantity, that $V(Z) = P(Z^2)\,-\,[P(Z)]^2$.



\section{Kullback's Symmetric Divergence Complex}

In this Section we shall learn how and why the structure of Kullback's symmetric divergence function naturally suggests that it be embellished by vectorial components in three independent dimensions.\\

\indent To begin, recall the details of the Kullback-Leibler {\it directed} divergence function $D(^.||^.)$ and its extropic complement $D^c(^.||^.)$.  These functions are also known as the relative entropy and extropy functions of $p$ with respect to $q$:
\begin{small}
\begin{eqnarray}
D(p||q) \ \ \ = \ \ \ \ \sum p \; log(\frac{p}{q}) \ \ \ \ \ \ &=& \ \sum p\; [log(p)-log(q)] \ > \ 0\; \nonumber {\rm ,\ and} \\ 
D^c(p||q) \ =  \ \sum (1-p) \; log(\frac{1-p}{1-q}) \ &=& \ \sum (1-p)\; [log(1-p) - log(1-q)] \ > \ 0 \ ; \nonumber 
\end{eqnarray}
\end{small}
\noindent and similarly for the ``reverse directions'',
\begin{small}
\begin{eqnarray}
D(q||p) \ \ \ = \ \ \ \ \ \sum q\; log(\frac{q}{p}) \ \ \ \ \ \ &=& \ \sum q\; [log(q)-log(p)] \ > \ 0 \nonumber \ {\rm,\ and} \\
D^c(q||p) \ =  \ \sum (1-q)\; log(\frac{1-q}{1-p}) \ &=& \ \sum (1-q)\; [log(1-q) - log(1-p)] \ > \ 0 \nonumber \ \ .
\end{eqnarray}
\end{small}
\noindent Those strict exceedance inequalities presume only that the pmvs $p,q$ are not identical:  $p\; \neq \;q$.  Kullback's (1959, p.6) symmetric divergence function $\mathscr{D}(p,q)$ is generated by summing the non-symmetric directed divergences: \  $\mathscr{D}(p,q) \ \equiv \ D(p||q)\;+\;D(q||p)$.  Similarly for the complementary symmetric divergence, \ $\mathscr{D}^c(p,q) \ \equiv \ D^c(p||q)\;+\;D^c(q||p)$.  The ``total symmetric divergence'' is denoted by  $\mathscr{D}_T(p,q) = \mathscr{D}(p,q) + \mathscr{D}^c(p,q)\; .$

\subsection{Expected comparative gains {\it sum} to the total divergence}

Now an intriguing idea derives from examining once again the expected comparative gain equations, this time algebraically in the context of the Total Log Scoring function, and recognising their relations to the Kullback-Leibler {\it directed divergence} functions.  We state this 
as a \\


\noindent {\bf Corollary 1 to Theorem 3.}\ \ \ In the context of the Total Log Scoring function $S_{TL}(X, p)$ introduced in equation (1), the previsions for the comparative gain reduce algebraically by direct substitutions to \vspace{.22cm}\\ $P_p[CG(X,p,q)] \, = \, \sum p \; log(\frac{p}{q}) \, + \,\vspace{.22cm} \sum (1-p) \; log(\frac{1-p}{1-q}) \, = \, D(p||q) + D^c(p||q) $ \;;\\
\noindent and similarly, \vspace{.22cm} \hspace{12.5cm}$(\ref{eq:prevscompgain}^*)$\\
 $P_q[CG(X,q,p)] \ = \ \sum q \; log(\frac{q}{p}) \ + \ \sum (1-q)\; log(\frac{1-q}{1-p}) \ = \ D(q||p) + D^c(q||p) $ \; .\\

\noindent Each pmv proponent's expected comparative gain  over the other's is equal to the relative entropy plus the relative extropy of the proponent's own pmv with respect to the pmv of the other.  This is the Kullback-Leibler directed distance, $D(^.||^.)$, plus its complementary dual, $D^c(^.||^.)$, both of which are known to be positive.  So both the proponents of $p$ and $q$ expect a positive comparative gain in total score relative to the other.  The positivity of each's expected gain over the other does not depend on which of them expects to receive a larger score personally.  Again, both of them expect to receive a larger score than the other (on account of the propriety of the scoring rule) but surely the size of only one of their own actual scores will be larger than that of the other.  In fact, their actual comparative gains are the negative equivalents of one another, according to their defining equation $(\ref{eq:compgain})$.\\


It is now worth noting that the {\it sum} of p's and q's previsions for their comparative gains 
can be seen in equations $(\ref{eq:prevscompgain})$ to equal the {\it total} symmetric divergence of the pmv's. \\

\noindent {\bf Corollary 2 to Theorem 3.} \ \ Again in the case of the Total Log Score,
\begin{small}
\begin{eqnarray}
P_p[CG(X,p,q)] + P_q[CG(X,q,p)] \ &=& \  D(p||q) + D^c(p||q) \;+\;D(q||p)+D^C(q||p) \nonumber \\ 
&=& \; \mathscr{D}(p||q) \; + \; \mathscr{D}^c(p||q) \nonumber \\
&\equiv& \; \mathscr{D}_T(p||q)\;. \nonumber
\end{eqnarray}
\end{small}
\noindent This is despite the fact that the sum of their {\it actual} comparative gains equals $0$.  Further investigation yields still another intriguing equivalence.

\subsection{Expected {\it Net Gains} sum to the same thing}

Intrigue arises from examining now the algebraic details of p's and q's expectations of each others' {\it net gains}, viz.,
\begin{small}
\begin{eqnarray}
P_p\{S(X,q) - P_q[S(X,q)]\} \ &=& \ \sum (p-q) \; log(q) \ + \ \sum (1-p-(1-q)) \; log(1-q) \nonumber \\
&=& \ \sum (p-q)\; log(\frac{q}{1-q}) \nonumber \ , \ {\rm \ \ and\ similarly}\hspace{2cm}\\
P_q\{S(X,p) - P_p[S(X,p)]\} \ &=& \ \sum (q-p)\; log(\frac{p}{1-p}) \nonumber \ \ .
\end{eqnarray}
\end{small}
\noindent Each of these is the difference between one's own
expection and the other's expectation of the other's log odds for
the occuring event.  \\

Now summing {\it these} expected net gains yields a companion to Corollary 2 of Theorem 3, evident from basic definitions:\\

\noindent {\bf Corollary 3 of Theorem 3.\ \ } Again using the Total Log Scoring rule, \vspace{.22cm}\\ \hspace*{3cm}$P_p[NG(X,q)] \ + \ P_q[NG(X,p)] \ = \ -\mathscr{D}_T(p||q)$\ .

\subsection{An alternative generator of  $\mathscr{D}_T(p||q)$.}

Perhaps surprisingly, at this late stage in the study of the divergence function, we can now recognise that still two more companion functions can be identified whose sums yield the same symmetric divergence as do the Kullback-Leibler functions $D(p||q)$ and $D(q||p)$. The awareness arises from Corollaries 2 and 3 of Theorem 3.  The first pair, also  directed functions, arises formally from switching the components of the finite discrete entropy function, $-H(p) = \sum p\;log\;p$ that get ``differenced'' in the production of $D(^.||^.)$.  Rather than thinking about $\sum p\;(log\;p - log\;q)$, we consider differencing the ``$p$'' instead of the ``$log\;p$'', raising consideration of $\sum (p - q)\;log\;p$.  We shall denote this alternative directed function by $\Delta(^.||^.)$, motivating it below:
\begin{eqnarray}
\Delta(p||q) \ \ &\equiv& \ \sum \; (p-q) \; log(p) \  \ \ , \ {\rm and} \label{eq:deltadef}\\
\Delta^c(p||q) \ &\equiv&  \ \sum (1-p-(1-q)) \; log(1-p) \ =  \ \sum (q-p) \; log(1-p) \ \ \ ; \nonumber 
\end{eqnarray}
\noindent and then similarly for the ``reverse directions'',
\begin{eqnarray}
\Delta(q||p) \ \ &=& \ \sum (q-p) \; log(q) \ \ \ , \ \ {\rm and} \hspace{6cm}\nonumber \\
\Delta^c(q||p) \ &=&  \ \sum (p-q) \; log(1-q) \ \ \ . \nonumber
\end{eqnarray}
\noindent Referring to our formulation of cross prevision assertions, it is apparent that $\Delta(p||q) \ =\ CH(q,p)-H(p)$, the amount by which $q$'s cross entropy for $p$ exceeds the entropy in $p$ itself.  Similarly,  $\Delta^c(q||p) \ = \ CJ(q,p)-J(p)$, the amount by which $q$'s cross extropy for $q$ exceeds $q$'s own extropy.\\

These considerations yield a surprising result: \ both the directed distance $D(p||q)$ and its reverse $D(q||p)$, and this new function $\Delta(p||q)$ and its reverse $\Delta(q||p)$, {\it sum} to the same thing, ...  Kullback's symmetric divergence, which we are denoting by $\mathscr{D}(p||q)$.  The same is true of the complementary sums, which we denote by $\mathscr{D}^c(p||q)$.  This result arises as  \vspace{.3cm}\\
\noindent {\bf Corollary 1 to Theorem 4.}\vspace{.3cm}\\
$\Delta(p||q) \; +\;  \Delta(q||p) \ = \ \sum (p-q) \; [log(p)-log(q)]\  = \vspace{.05cm} \\ 
\hspace*{3cm} = \ D(p||q)\;  + \; D(q||p) = \ \mathscr{D}(p,q)$, \vspace{.22cm}\\
\noindent and \vspace{.22cm}\\
\noindent $\Delta^c(p||q)  +  \Delta^c(q||p) \ = \ \sum (q-p) [log(1-p) - log(1-q)]  \ = \vspace{.05cm} \\
\hspace*{3.1cm} = \ D^c(p||q) +  D^c(q||p) \ = \ \mathscr{D}^c(p,q)$.\\

\noindent Here the sums of the $\Delta(^.||^.)$ functions apply the total log score to the Net Gains, while the $D(^.||^.)$ functions arise from  its application to the Comparative Gains.


\subsection{Identifying still a third function that yields $\mathscr{D}_T(p,q)$}

Expanding the identical sums $D(p||q)+D(q||p)$ and $\Delta(p||q)+\Delta(q||p)$ which constitute $\mathscr{D}(p,q)$ into their component pieces, and then exhibiting them across a line makes it apparent that each of the  functions, $D(p||q)$ and  $\Delta(p||q)$  arises by selecting two of the four summands to define the function, and thus the remaining two to define their ``reverses'',  $D(q||p)$ and  $\Delta(q||p)$.\vspace{.12cm}\\

\includegraphics[width=0.75\textwidth]{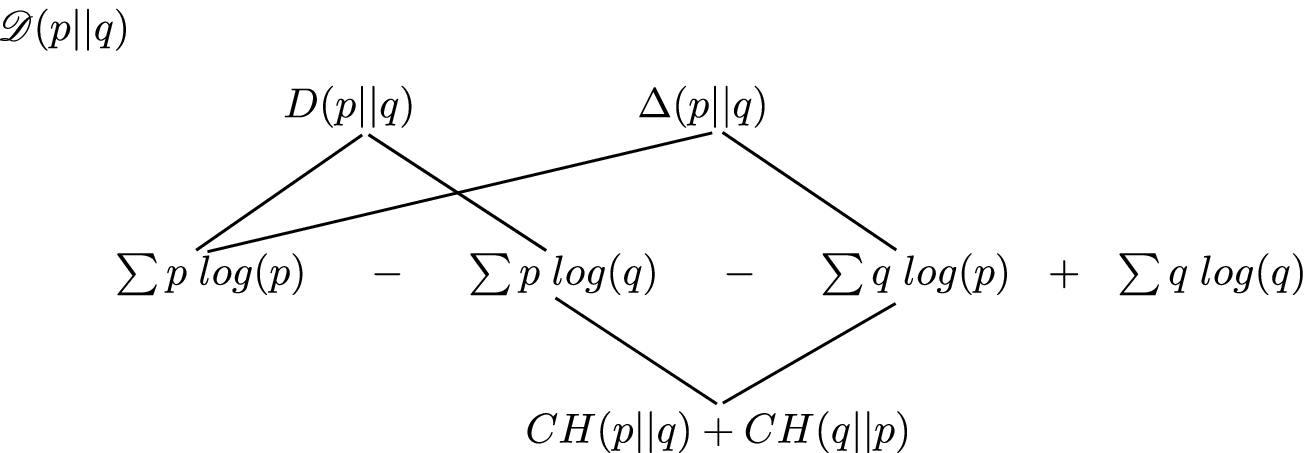}\vspace*{.38cm} \\
\noindent  $D(p||q)$ is the sum of the first two pieces of the first line, leaving $D(q||p)$ composed of the second two;  while $\Delta(p||q)$ is the sum of the first and third piece, leaving $\Delta(q||p)$ to be composed of the second and fourth piece. \\

Now seeing these sums and these function pairs in this way suggests that a third (and the accompanying sixth) choice of two of the four summands may be worthy of study as well.  Consider using the interior two summands (the second and third pieces of the first line) to constitute still another form of function that would sum (along with the remaining pair) to yield the symmetric divergence $\mathscr{D}(p||q)$.  This choice would be the sum of the ``cross entropy'' functions, \ $CH(p||q)  +  CH(q||p)$.  One difference of this function from $D(p||q)$ and $\Delta(p||q)$ is that this sum of cross entropies is already symmetric on its own, whereas the latter two functions require summation with their companions $D(q||p)$ and $\Delta(q||p)$ to construct a symmetric function. A second difference of the sum $CH(p||q)+CH(q||p)$ from the $D(p||q)$ and $\Delta(p||q)$ functions appears in its location for self-divergence.  Whereas the values $D(p||p) \ = \ \Delta(p||p) \ = \ 0$, the value of $CH(p||p)+CH(q||q) \ = \ H(p)+H(q)$, the summed entropies in $p$ and in $q$.  However, defining the symmetric divergence function more generally than heretofore so to eliminate this location shift of ``self-distance'' yields a pleasing symmetrical result.\vspace{.4cm}\\
\noindent{\bf Theorem 5.} \ \ Based on a generalised definition of symmetric divergence, so to account for the location of ``self-divergence'',
\begin{eqnarray}
\mathscr{D}(p,q) \ &\equiv& \ D(p||q)\;+\;D(q||p)\;-\;D(p||p)\;-\;D(q||q) \nonumber \\
 &=& \ \Delta(p||q)\;+\;\Delta(q||p)\;-\;\Delta(p||p)\;-\;\Delta(q||q)  \label{eq:gensymndiverg}\\
  &=& \ CH(p||q)\;+\;CH(q||p)\;-\;CH(p||p)\;-\;CH(q||q) \nonumber \\
  &\geq & \ 0 \;.\nonumber
\end{eqnarray} 
\noindent Of course those self-divergence values in the first two lines equal $0$, whereas in the third line $CH(p||p) = H(p)$ and $CH(q||q) = H(q)$.\\

The same procedures of deconstruction and reconstruction can be performed with the complementary functions $(D^c(^.||^.)$ and $ \Delta^c(^.||^.))$ in a second line: \ $D^c(p||q)$ is the sum of the first two pieces of this second line, while $\Delta^c(p||q)$ is the sum of the first and third pieces. \\

\includegraphics[width=0.95\textwidth]{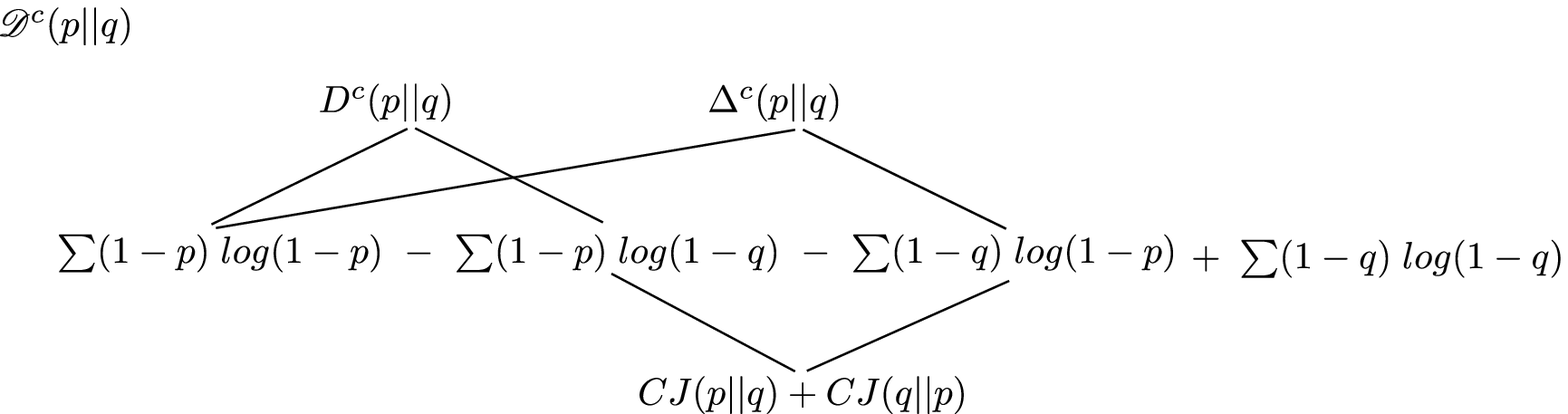}\vspace*{.38cm} 

We propose that the Kullback directed distance function needs to be appended by all three of its three generators if it is to characterise the information content of 
the statistical problem the forecasters address.  Together they constitute a complex. \\

One might be puzzled. ``Who cares about all these generators of the same symmetric distance?  The symmetric distance is what it is.''  The answer derives from an identifiable linear relation between the four functions 
$\mathscr{D}_T(p,q),\ D(p||q) + D^c(p||q), \ \Delta(p||q) + \Delta^c(p||q)$, $CH(p||q) + CJ^c(p||q)$ \ and the ``four fundamental previsions'' we introduced in equation (\ref{eq:4fundprevs}). We address this now.  

\subsection{The isomorphism of the Kullback complex of generators with \\
the four fundamental previsions}

Each of these eight component functions in a Kullback complex is a different linear combination of the various entropies, extropies, cross entropies and cross extropies that have been central to our scoring analysis.  These combinations happen to be ordered in such a way that the two spaces spanned by the 4-dimensional vector functions are isomorphic, related by the linear equations 
 \begin{center}
$\left[ \begin{array}{c}  \mathscr{D}(p||q)+\mathscr{D}^c(p||q) \\ D(p||q)+D^c(p||q) \\ \Delta(p||q)+\Delta^c(p||q) \\
C_H(p||q)+C_J(p||q)\\
\end{array} \right] = \begin{bmatrix} 1&-1&1&-1 \\ 1&-1&0&\ \ 0 \\1&\ \ 0&0&-1 \\0&-1&0&-1 \end{bmatrix} \times \left[ \begin{array}{c} P_p[S(X,p)]\\ P_p[S(X,q)]\\ P_q[S(X,q)] \\P_q[S(X,p)] \label{eq:isomorphism}  \end{array} \right] \ \ (8)$ \\
\end{center}
\noindent which is easily confirmed, where $C_H(p||q)\equiv CH(p||q)+CH(q||p)$ and 
$C_J(p||q)\equiv CJ(p||q)+CJ(q||p)$.  The inverse transformation is
\begin{center}
$\left[ \begin{array}{c} P_p[S(X,p)]\\ P_p[S(X,q)]\\ P_q[S(X,q)] \\P_q[S(X,p)] \end{array} \right] = \begin{bmatrix} 0&\ \ .5&\ \ .5&-.5\\0&-.5&\ \ .5&-.5\\1&-.5&-.5&-.5\\0&\ \ .5&-.5&-.5\end{bmatrix} \times \left[ \begin{array}{c} \mathscr{D}(p||q)+\mathscr{D}^c(p||q) \\ D(p||q)+D^c(p||q) \\ \Delta(p||q)+\Delta^c(p||q) \\ C_H(p||q)+C_J(p||q)\\
\end{array} \ \  \right] . \  (9)$\vspace{.3cm}\\
\end{center}

What this relation tells us is that each of the three distance function generating pairs measures a distinct information source for understanding of the full information content of Kullback's symmetric distance. The symmetric distance measure  is incomplete on its own: \ it needs to be supplemented by three more components if it is to represent the information content of the pmv assertions p and q.  In all, it is identifiable better as a {\it vector} which we might refer to as the Kullback complex ... $[\mathscr{D}(p,q), \ D(p||q), \Delta(p||q), CH(p||q)+CH(q||p)]$, rather than merely by the symmetric function value $\mathscr{D}(p,q)$ alone.\\

One further recognition is in order.  We have couched this discussion in the context of the total logarithmic score $(\ref{eq:totlogscore})$ and the Kullback total symmetric distance function $\mathscr{D}_T(p||q)$.  It is worth remarking that the total log score can be considered to be the sum of two separate components: the log score itself, $\sum E_i\;log\:p_i$, and the complementary log score, $\sum \tilde{E_i}\;log\:(1-p_i)$.  Actually, both of these components of the total log score are proper scores on their own, as is {\it any} linear combination of them.  The two component scores assess different aspects of the asserted pmv's $p$ and $q$.  One piece focuses more directly on the center of the mass function, while the other focuses more on the tails.  Because these aspects of a pmv may concern different people in different ways with respect to utility, it is usually best to keep track of the two components of the score separately.  If the scoring function $S_{TL}(X,{\bf p}_N)$ were considered to be the sum of the two separate scores, say, $S_{log}(X,{\bf p}_N)$ + $S_{complog}(X,{\bf p}_N)$, the isomorphism equations $(9)$ hold separately for each component of the total log score and their previsions too.  Tracking these components separately, the total symmetric Kullback complex could be considered to be an 8-dimensional space of measurement functions. \\

We shall now address a final technical note concerning the relation of these newly recognised dimensions of the information complex to the theory of Bregman divergences.

\section{$\Delta(^.||^.)$ and $HC(^.||^.)$ are {\it not} Bregman divergences}
Information theorists are comfortable by now with the fact the Kullback's symmetric divergence function is the sum of two ``directed divergences'', each summand of which is a Bregman divergence.  That is, $\mathscr{D}(p||q) \; = \; D(p||q) + D(q||p)$, and $D(p||q)$ is a Bregman divergence with respect to the Bregman function $\Phi(p) = \sum p\;log\;p$.  The same is true of the complementary Kullback symmetric divergence:  \ $\mathscr{D}^c(p||q) \; = \; D^c(p||q) + D^c(q||p)$, where $D^c(p||q)$ is the Bregman divergence associated with the Bregman function $\Phi^c(p) = \sum (1-p)\;log\;(1-p)$.  We shall remind you of the definition of a Bregman divergence shortly.\\

It may come as a surprise then that 
{\it neither} of the companion generators of Kullback's symmetric function that we have found constitutes a Bregman divergence. These are the functions we have denoted by $\Delta(p||q) \equiv \sum (p-q)\;log\;p$,\ and $CH(p||q) \equiv \sum p\;log\;q$.  In brief, a function $d_{\Phi}(p,q)$ is a Bregman divergence with respect to a convex, differentiable function $\Phi({\bf p}_N)$ if $d_{\Phi}(p,q)\;=\;\sum\;[\phi(p)-\phi(q) - (p-q)\phi'(q)]$.  To be precise, this describes only the ``separable'' context in which $\Phi({\bf p}_N)$ has the form $\sum \phi(p_i)$, where $\phi(^.)$ is a function of the single component variables of the vector $p$.  However, this {\it is} the limited context that is relevant to our discussion here. \\

Using some simple algebraic adjustments then, $\Delta(p||q)$ and $CH(p||q)$ would also be identifiable as Bregman functions if and only if there were solutions $\phi(^.)$ to the following two problems of differential equations:\\

\noindent 1.\ Find a convex differentiable function $\phi(^.)$, with a domain that  includes both $p$ and $q$ in the interval $(0,1)$, and for which
\begin{center}
 \ $\sum \;(p-q)\; log \; p \ = \ \sum\; [\phi(p)-\phi(q) - (p-q)\phi'(q)]$
\end{center}
\noindent or equivalently for which, \ \ $\phi'(q) \ = \ \frac{\phi(p) - \phi(q)}{p-q} \; - \; log\;p$\; \ for any $p$ and $q \;\in\;(0,1)$\ ; \vspace{.22cm}\ \ and \\

\noindent 2.\ \ Find another such convex differentiable function $\gamma(^.)$ for which
\begin{center}
\ $\sum\; p\; log \; q \ = \ \sum\; [\gamma(p)-\gamma(q) - (p-q)\gamma\:'(q)]$
\end{center}
\noindent or equivalently for which, \ \ $\gamma\:'(q) \ = \ \frac{\gamma(p) - \gamma(q)}{p-q} \; - \; \frac{p\;log\;q}{p-q}$\; \ for any $p$ and $q \;\in\;(0,1)$\ . \\

\noindent If such functions could be found, then they would be the Bregman functions that specify $\Delta(p||q)$ and $CH(p||q)$ as their Bregman divergences, respectively.\\

However, neither of these problems does have solutions!  This is revealed by examining the limits as $p \rightarrow q$ of these two conditioning equations on the sought-for functions. In this limiting context, the first problem would require that $\phi'(q)\;=\;\phi'(q) + log\;q$; \; while the second would require another nonsense ... that $\phi'(q)\;=\;\phi'(q) + \infty$.  On the contrary, if the left-hand-side of the conditions of problems 1 or 2 were replaced by Kullback-Leibler's function $D(p||q)$, then this limiting condition would require only that $\phi'(q)\;=\;\phi'(q) + lim_{p \rightarrow q}\; \frac{p\;log\;p - p\;log\;q}{p-q}\;-\;1$\:, a condition that is assured by L'Hopital's rule. This feature is what distinguishes the Kullback-Leibler measure as a Bregman function.\\

The upshot of this understanding is that the dissection of of $\mathscr{D}(p,q)$ via each of the three directed functions $D(^.||^.), \Delta(^.||^.), $ and $CH(^.||^.)$ illuminates a distinct dimension of the measure we know as Kullback's total symmetric distance.  We might well call the complete vector of these measures a ``Kullback complex.''  All arguments regarding the accompanying dimensions of $\mathscr{D}_T(p||q)$ pertain individually to its direct and complementary components, $\mathscr{D}(p||q)$ and $\mathscr{D}^c(p||q)$, as well.\\

A concluding query arises then from this analysis.  What is the real motivating feature of Bregman divergences that makes them so widely well regarded?  We have found out that two other equally fundamental measures generate the same symmetric divergence as does the Kullback divergence, but these newly recognised generating functions are {\it not} Bregman divergences.  Evidently they do provide substantive content to a complex of measures which join with the symmetric divergence itself to summarize information.

\section{The subjectivist take on the Kullback functions}

Kullback's classic work (1959) on {\it Information Theory and Statistics} developed out of quite a  different statistical imagination than the subjectivist construction we have been developing here.  He began (pp. 4-5) by positing two hypothesized random generating functions of an observation, $X$, described by the densities $f_1(x)$ and $f_2(x)$.  Ignoring here some measure theoretic refinements, he then identified the measure $\int f_1(x)\;log(\frac{f_1(x)}{f_2(x)})\;dx \ \equiv \ I(1:2)$ as ``the information in $X=x$ for discriminating in favour of the generator $f_1(x)$ as opposed to $f_2(x)$.''  This choice was motivated by the measure's role in distinguishing the expected posterior log odds from the prior log odds in favour of $H_1$ relative to $H_2$.  It was this nonsymmetric ``information function'' $I(1:2)$ which he joined with its reverse function $I(2:1)$ to yield his symmetric divergence function $J(1,2) = I(1:2) + I(2:1)$.\\

Indeed, this framework currently still provides the standard objectivist context for most all contemporary work in information statistics.  The subjectivist context in which we honour the pathbreaking insight of Bruno de Finetti construes all statistical problems quite differently.\\

In the subjectivist context of the work we have been addressing here, the notation for Kullback's symmetric divergence function $J(1,2)$ has been replaced notationally by $\mathscr{D}(p,q)$, and is used for another purpose. The functions $f_1(x)$ and $f_2(x)$ are replaced by our vectors $p$ and $q$.  They are understood to represent not ``unobservable potential random generators of $X$'' that need to be ``tested'', but as pmv vectors of  personally asserted probabilities for the possible values of $X$ by two different people based on their differing information sources and assessments.
Kullback's non-symmetric information function $I(1:2)$ has been denoted herein as  $D(p||q)$.  We now find that two other  functions, $\Delta(p||q)$ and $CH(p||q)$, can be used to generate the same symmetric divergence function $\mathscr{D}(p||q)$ as does $D(p||q)$ when it is defined more completely, and we know why.  The Kullback complex is merely a translation of the information measures encoded by our four fundamental expectation equations into  different ranges and scales.  The information components $D(^.||^.)$, $\Delta(^.||^.)$, and $CH(^.||^.)$ are linearly distinct from $\mathscr{D}(^.||^.)$ itself.


\section{On statistical application}

We have appreciated the opportunity to present this analysis as a tribute to the marvelous insights provided by the work of Bruno de Finetti.  It is only the limitation of space that preclude a complete exhibition of a substantive application of its implications here.  Two remarks are in order.

\subsection{A preliminary report}
The authors have been engaged for several years in various applications of the comparative evaluation of empirical forecasting distributions using proper scoring rules.  In so doing we have used the the traditional style of awarding forecasters the proper score appropriate to their own forecasts.  However, we had been bothered by the arbitrary aspects involved in that method that we have mentioned in Section 2.2.
In collaboration with our colleague Gianna Agr\`o, we are currently engaged in a reassessment of a proper scoring analysis of alternative forecast distributions for daily stock returns that exhibit strikingly different structures in the tail properties of the forecasting distributions. Our original analysis conducted in the traditional style was reported in  Agr\`o et al. (2010). Though our work on the reassessment is yet  incomplete, we are already aware that the implications of the Pareto exchange methodology make major differences in the assessment of the two forecast distributions.
An analysis of the complete results of this research will soon be accessible.  It is currently only in working form in Sanfilippo et al. (2017).  \\

To provide a taste of the strength of the difference the Pareto exchange makes in applications, we present here in Figure 1 a side-by-side comparison of the cumulative Pareto-exchange scores and the simple cumulative personal scores for two mass functions.  Without further description here, they are identified only as Forecast density $f(x_{t+1}|{\bf x}_t)$ and Forecast density $g(x_{t+1}|{\bf x}_t)$.  The differences in the historical structures of the scoring behaviour are startling.  Rather than suggesting a regular improvement in the forecast of ``g'' upon that of ``f'' as seen in Figure~1b according to Personal Direct Log Scoring, the Cumulative Pareto Exchange Log Score displays changing periods of dominance by each of the alternative forecasts in Figure 1a.  
Moreover, these are left in a cumulative state that favours the Forecast ``f'' rather than "g".   
A complete specification of precise details must be delayed to the full report.
\begin{figure}[ht]
    \captionof{figure}{ {\bf Comparative results of Pareto exchange Scores and Direct Scores}}
        \centering
        \begin{subfigure}{0.486\textwidth}
                \centering
                \includegraphics[width=1\linewidth]{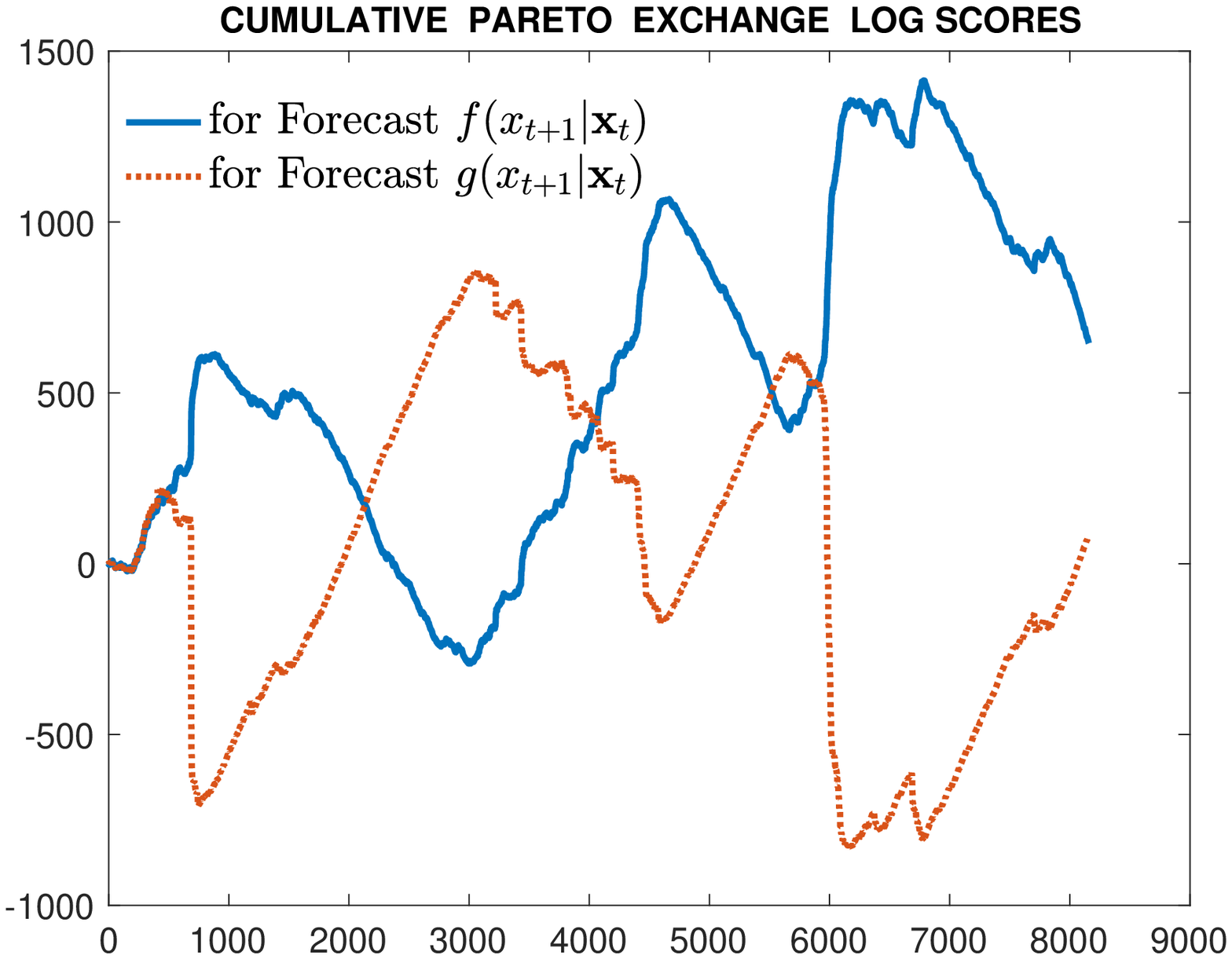}
                \caption{The Pareto exchange of normed Log
                        Scores identifies periodic shifts in the
                        assessed quality of forecast performances
                        by ``f'' and ``g''. }
                \label{fig:1a}
        \end{subfigure}
        \hspace{0.22cm}
        \begin{subfigure}{0.48\textwidth}
                \centering
                \includegraphics[width=0.98\linewidth]{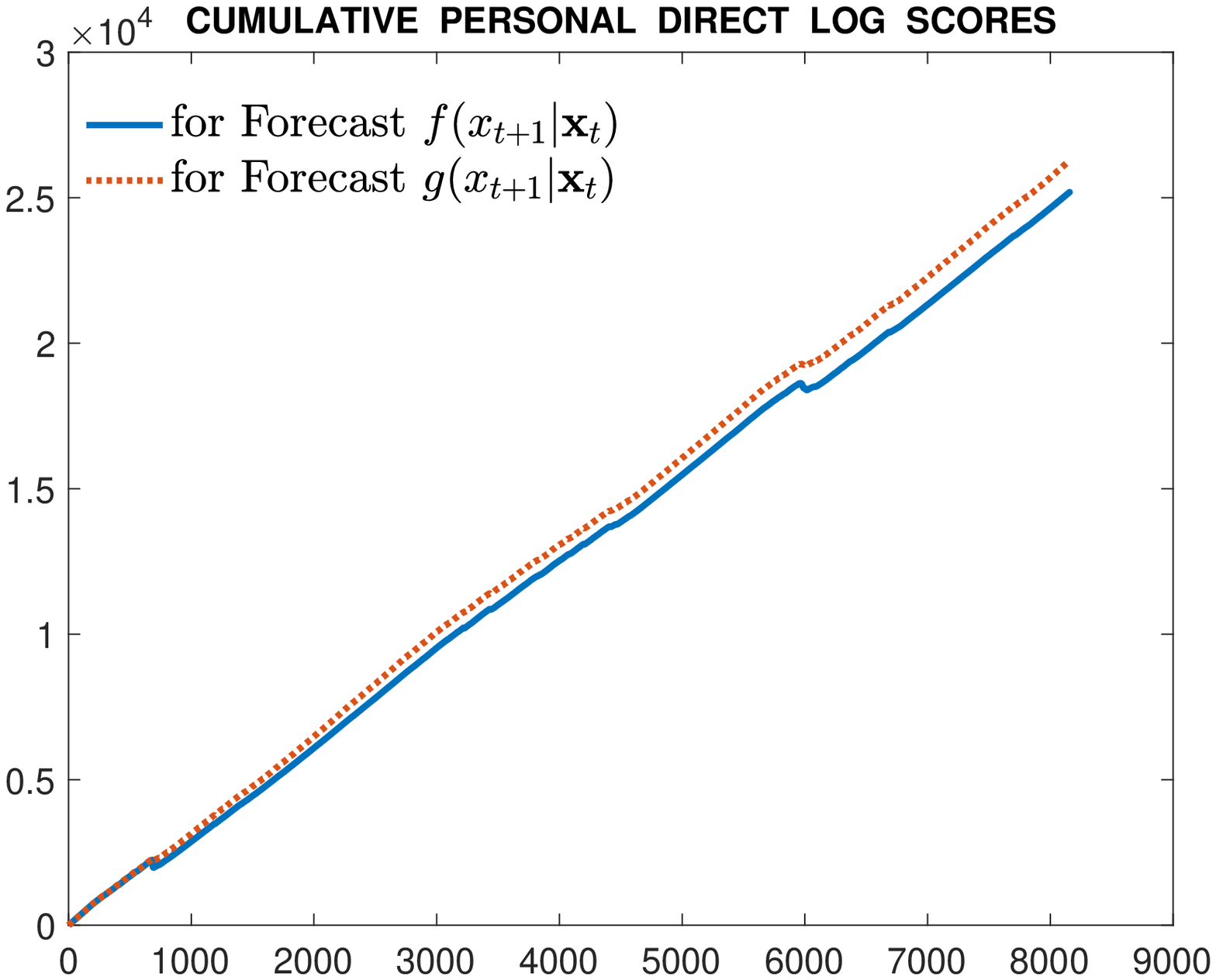}
                \caption{The direct accumulation of Log Scores 
                      suggests
                     regularly improved quality of forecast ``g''
                        relative to that of ``f''.\\}
                \label{fig:1b}
        \end{subfigure}
        \label{fig1}
\end{figure}

\subsection{The subjectivist alternative to  hypothesis testing}

In our introductory remarks we mentioned that from de Finetti's subjectivist perspective, the time honoured practice of statistical hypothesis testing 
needs to be replaced by methods based on proper scoring rules.  Without presenting an exemplary  analysis, we can sketch here the structure of such a procedure based on the Pareto optimal exchange of proper scores.\\

Consider the context of the simplest statistical problem for which the Neyman-Pearson Lemma applies.  A sequence of quantities presumed to be independent emanations of a distribution with mean $\mu$ are to be  observed, and it is proposed to test the null hypothesis $H_0$ that $\mu = \mu_0$ against an alternative $H_A$ that $\mu = \mu_A$.  Despite early recognition in Lehmann (1959, p 61; p 69 in Second Edition) that the choice of an appropriate ``significance level'' required for performing such a test needs be based on a decision maker's relative utility for the correct and erroneous choice of the two hypotheses, statistical practice has drifted toward the understanding of arbitrary significance levels $.05, .01$, and $.001$ as summary guides to an appropriate decision to ``reject'' or to ``accept'' the null hypothesis.  More egregiously, the summary has been long extended to the computation and report of a ``p-value''.\\

Bruno de Finetti's subjectivist outlook would characterise this situation quite differently.  The sequence of observations are recognised as meriting a judgment of exchangeability, represented by a distribution which is necessarily then  mixture-hypergeometric.  The relative forecast quality of two such distributions could then be evaluated by methods exposed here.  One would be based upon a mixture favouring the location $\mu_0$, and the other favouring $\mu_A$.  The accumulating sequential scores of their forecasts would involve a Pareto exchange of their  proper scores.
An exposition of a simple computational example in this standard context would be appropriate to introductory level  texts on applied statistics.

\section{Concluding remark}

A little more than 50 years have passed since Buno de Finetti voiced his ``partially baked'' musing in I. J. Good's anthology:  ``Does it make sense to speak of 'good probability appraisers'?''  He had wondered about this question at that time because of the subjectivist attitude, which we support, that no probability assertion can be determined to be ``wrong'' on the basis of observational evidence.  For the probability assertion merely constitutes an honest expression of its promoter's uncertain knowledge about the value of the unknown quantity.  There is nothing wrong with being uncertain.  That is our common state.   Nonetheless, the analysis we have provided in this tribute essay supports an answer of ``Yes!'' to de Finetti's topical question. While details of our analysis pertain to the context of the total log score, our identification of the Pareto exchange is applicable to any proper scoring function.   Proper scoring rules {\it can} be used 
meaningfully to aid in the evaluation of the information quality in personal 
probability assertions as well as in elicitation of personal probabilities themselves.  However, the harvest of this usefulness requires both a switch in the scores that have commonly been awarded to each of the mass function proponents and an expansion of our understanding of the symmetric Kullback divergence measure.  This  has long been considered to be the cornerstone standard valuations of statistical information. 

\section*{Acknowledgements}
Thanks to Rua Murray for a suggestion regarding the Bregman derivatives, and to Gianna Agr\`o for her permission to print here the preliminary Figure 1 which represents the result of our joint work in progress.  Thanks to Anuj Misra for producing the visual schema representing the construction of Kullback's symmetric divergences $\mathscr{D}(p,q)$ and $\mathscr{D}^c(p,q)$.

\section*{References}
\noindent{Agr\`o, G., Lad, F., and Sanfilippo, G.} (2010) Sequentially forecasting economic indices using mixture linear combinations of EP distributions, {\it Journal of Data Science}, {\bf 8}, 101-126.

\vspace*{.22cm}
\noindent{Dawid, A.P. and Musio, M.} (2015) Bayesian model selection based on proper scoring rules, {\it Bayesian Analysis}, {\bf 10}, 479-499.
\newpage
\noindent{de Finetti, B.} (1962) Does it make sense to speak of ``good probability appraisers''?,   in {\it The Scientist Speculates:  an anthology of partly-baked ideas}, I. J. Good (ed.), London: William Heinemann, 357-364.

\vspace*{.22cm}
\noindent{de Finetti, B.} (1967) Quelque conventions qui semblent utiles,  {\it Revue Romaines des Mathe\'matiques Pures et Applique\'s}, {\bf 12}, 1227-1233.  L.J. Savage (tr.) A useful notation, in {\it Probability, Statistics and Induction:  the art of guessing}, New York: John Wiley (1972), pp. xvii-xxiv.

\vspace*{.22cm}
\noindent{Gneiting, T. and Raftery, A.} (2007) Strictly proper scoring rules, prediction, and estimation, {\it JASA}, {\bf 102}, 359-378.

\vspace*{.22cm}
\noindent{Good, I.J..} (1962) {\it The Scientist Speculates: an anthology of partly-baked ideas}, London: \; William Heinemann.  This includes the article of de Finetti (1962).

\vspace*{.22cm}
\noindent{Kullback, S.} (1959) {\it Information Theory and Statistics}, New York: John Wiley.

\vspace*{.22cm}
\noindent{Kullback, S. and Leibler, R.A.} (1951) On information and sufficiency, {\it Annals of Mathematical Statistics}, {\bf 22}, 79-86.

\vspace*{.22cm}
\noindent{Lad, F.} (1996) {\it Operational Subjective Statistical Methods: a mathematical, philosophical, and historical introduction}, New York: John Wiley.

\vspace*{.22cm}
\noindent  {Lad, F., Sanfilippo, G., and Agr\`o, G.} (2012)  Completing the
logarithmic scoring rule for assessing  probability distributions, ISBrA Conference, AIP Conf. Proc. 1490, {\bf
13}, 13-30;  View online: 
http://dx.doi.org/10.1063/1.4759585 

\vspace*{.22cm}
\noindent{Lad, F., Sanfilippo, G., and Agr\`o, G.} (2015) Extropy: complementary dual of Entropy, {\it Statistical Science}, {\bf 30}, 40-58.

\vspace*{.22cm}
\noindent{Lehmann, E.L.} (1959) {\it Testing Statistical Hypotheses}, New York: John Wiley.

\vspace*{.22cm}
\noindent{Matheson, J. and Winkler, R.} (1976) Scoring rules for continuous probability distributions, {\it Management Science}, {\bf 22}, 1087-1096.

\vspace*{.22cm}
\noindent{Sanfilippo, G., Lad, F., and Agr\`{o}, G.} (2017)  Scoring alternative forecast distributions using Pareto exchanges of proper scores. Working paper.

\vspace*{.22cm}
\noindent{Savage, L.J.} (1971) Elicitation of personal probabilities and expectations, {\it JASA}, {\bf 66}, 778-801.

\end{document}